\documentclass[10pt,letterpaper]{article}

\usepackage[top=0.85in,left=2.75in,footskip=0.75in]{geometry}
\usepackage {indentfirst} 
\usepackage{amsmath,amssymb}

\usepackage{changepage}

\usepackage{textcomp,marvosym}

\usepackage{cite}

\usepackage{nameref}

\usepackage[right]{lineno}

\usepackage[nopatch=eqnum]{microtype}
\DisableLigatures[f]{encoding = *, family = * }

\usepackage[table]{xcolor}

\usepackage{array}
\usepackage{amsthm} 
\usepackage{comment}

\newcolumntype{+}{!{\vrule width 2pt}}

\newlength\savedwidth

\raggedright
\usepackage[aboveskip=1pt,labelfont=bf,labelsep=period,justification=centering,singlelinecheck=off]{caption}

\makeatletter
\renewcommand{\@biblabel}[1]{\quad#1.}
\makeatother
\usepackage{lastpage,fancyhdr,graphicx}
\usepackage{epstopdf}
\usepackage{float}
\usepackage{amsmath}
\usepackage{amssymb}
\usepackage[figuresright]{rotating}
\usepackage{hyperref}
\usepackage{graphicx}
\usepackage{subcaption}
\usepackage{booktabs}
\usepackage{tabularx}
\usepackage{comment}
\usepackage{ragged2e}
\usepackage{authblk}
\usepackage[none]{hyphenat}
\fancyheadoffset[L]{2.25in}
\fancyfootoffset[L]{2.25in}
\newcommand{\enabstractname}{Abstract}
\newenvironment{enabstract}{
\par\noindent\mbox{}\hfill{\bfseries \flushleft\enabstractname\hfill\mbox{}\par
\vskip 2.5ex}{\par \vskip 2.5ex}}

\title{Application of Superconducting Technology in the Electricity Industry: A Game-Theoretic Analysis of Government Subsidy Policies and Power Company Equipment Upgrade Decisions}
        \author[1]{Mingyang Li}
        \author[1]{Maoqin Yuan\thanks{Corresponding author: mqyuan@cupk.edu.cn}}
        \author[2]{Han Pengsihua}
        \author[3]{Yuan Yuan}
        \author[1]{Zejun Wang\thanks{Corresponding author: wangzj@cupk.edu.cn}}
        \affil[1]{School of Science and Arts, China University of Petroleum-Beijing at Karamay, Karamay 834000,Xinjiang, China.}
        \affil[2]{School of Business Administration, China University of Petroleum-Beijing at Karamay, Karamay 834000,Xinjiang, China.}
        \affil[3]{State key Laboratory of Petroleum Resources and Prospecting, China University of Petroleum, Beijing 102249, China.}
        \date{}
        
\begin{document}
	       \maketitle
	 \justify

        \begin{enabstract}
            This study investigates the potential impact of "LK-99," a novel material developed by a Korean research team, on the power equipment industry. Using evolutionary game theory, the interactions between governmental subsidies and technology adoption by power companies are modeled. A key innovation of this research is the introduction of sensitivity analyses concerning time delays and initial subsidy amounts, which significantly influence the strategic decisions of both government and corporate entities. The findings indicate that these factors are critical in determining the rate of technology adoption and the efficiency of the market as a whole. Due to existing data limitations, the study offers a broad overview of likely trends and recommends the inclusion of real-world data for more precise modeling once the material demonstrates room-temperature superconducting characteristics. The research contributes foundational insights valuable for future policy design and has significant implications for advancing the understanding of technology adoption and market dynamics.
        \end{enabstract}

	\section{Introduction}

	The energy-environment issue is a significant challenge in the 21st century, with electrical power loss in electricity companies being a particularly salient problem \cite{-1}. Power loss typically ranges from 5\% to 10\% in alternating current transmission lines and can exceed 15\% in some cases \cite{0}. As the largest electricity consumer in the world, China faces an even more challenging power loss issue in its electricity companies \cite{tongji}. Despite using ultra-high voltage transmission technology to minimize losses, China's annual energy loss in transmission lines is several hundred billion kWh. Therefore, technical innovations, such as superconductive technology, are crucial in the electricity industry \cite{2}.
	
	\begin{figure}[H]
		\centering
		\includegraphics[width=0.7\linewidth]{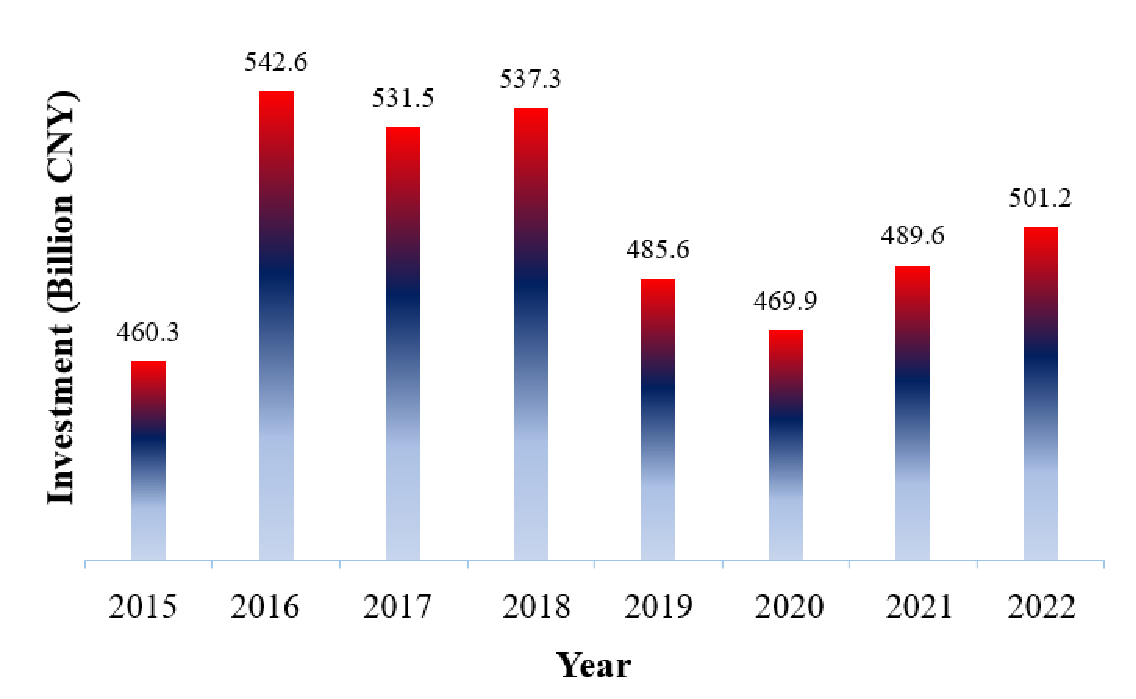}
		\caption{China power industry statistical data. Note: Ref \cite{tu2}}
		\label{0161}
	\end{figure}

	Superconductivity technology has significant applications in fields such as energy, manufacturing, and information and communication \cite{1}. It dramatically improves the efficiency of electrical energy transmission by allowing current to flow with zero resistance in materials under low-temperature conditions \cite{3}. Superconducting feeder cables have been proven capable of being introduced into railway feeder systems, reducing power transmission losses \cite{niu}.
	In July 2023, a research team from South Korea announced a breakthrough in room-temperature superconductivity, naming the material capable of room-temperature superconductivity "LK-99". This discovery garnered widespread attention globally\cite{lk}. Although the current validation results of this material do not fully meet expectations, its superconductivity and partial anti-magnetism exhibited under low-temperature conditions still hold significant research value and potential for application. The use of superconductivity technology in the power sector can effectively reduce power loss, thereby lowering the risk of overload, improving the quality of electrical energy, and enhancing grid stability. Superconductivity technology has brought potential impacts on the power grid, and it is expected to play a crucial role in the future smart grids \cite{4}.
	
	\begin{figure}[H]
		\centering
		\begin{minipage}[t]{0.45\textwidth}
			\centering
			\includegraphics[width=\linewidth]{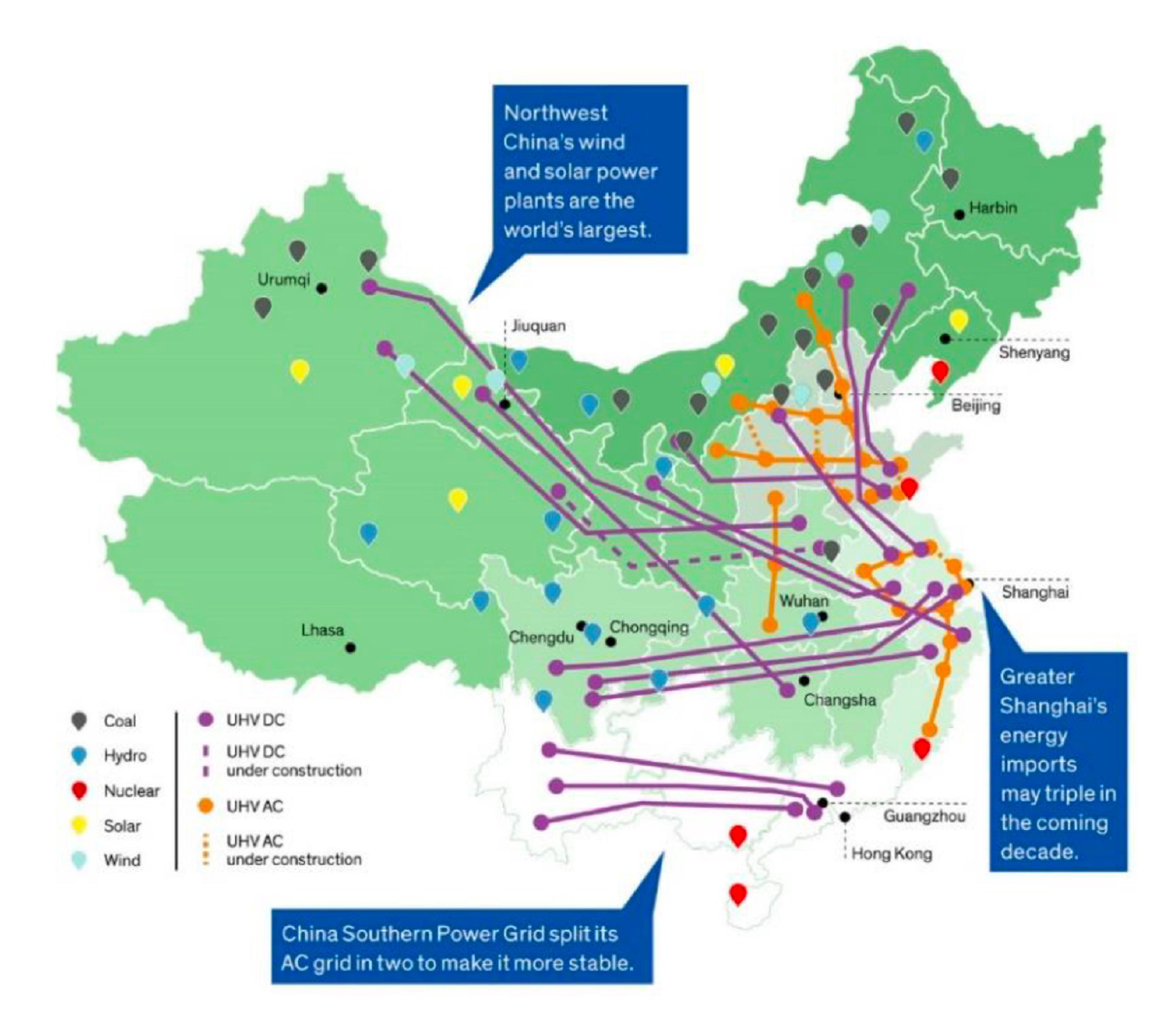}
			\caption{China’s hybrid AC-DC grids. Note: Ref \cite{31}}
			\label{fig:screenshot010}
		\end{minipage}%
		\hfill
		\begin{minipage}[t]{0.48\textwidth}
			\centering
			\includegraphics[width=\linewidth]{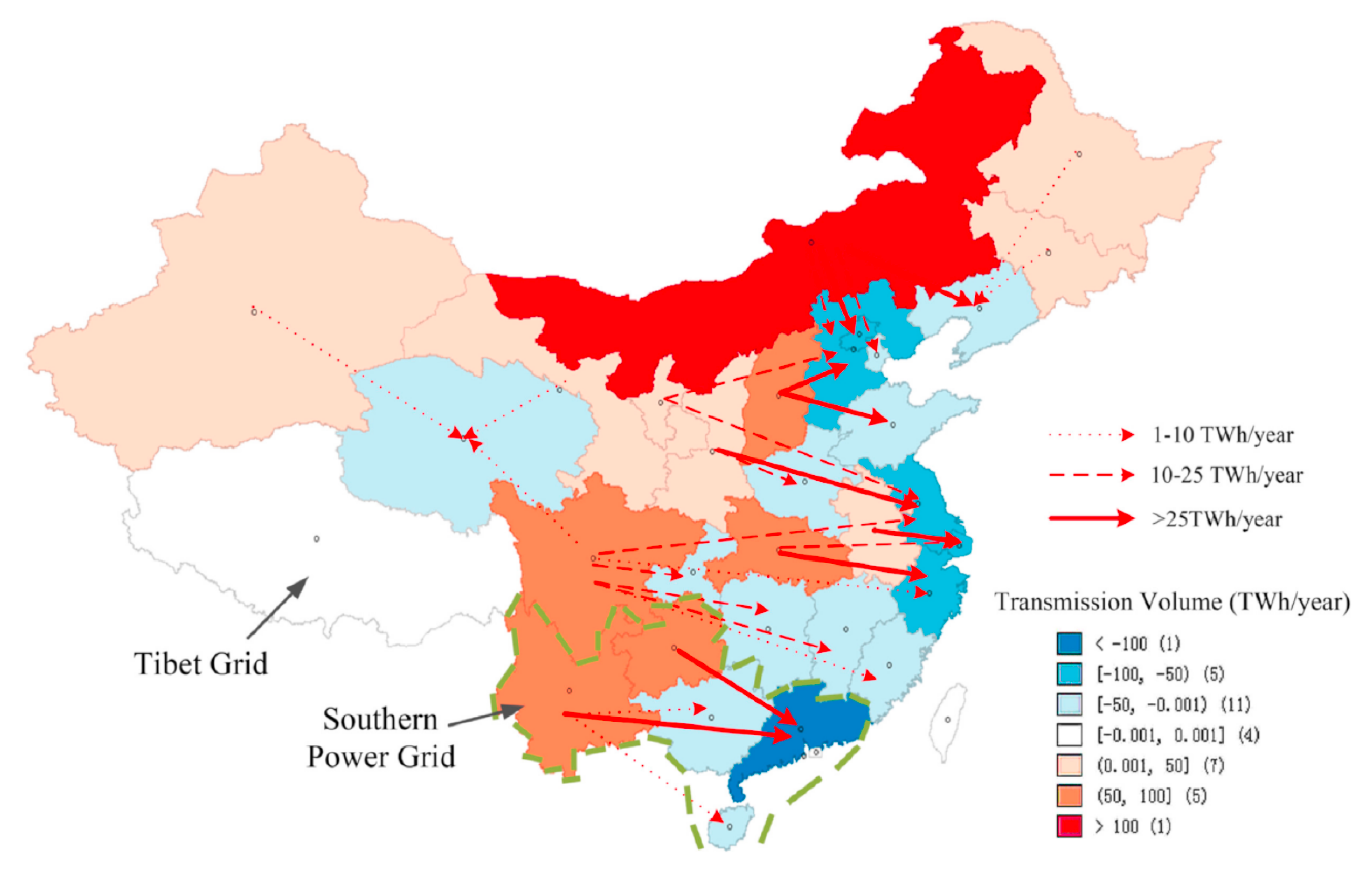}
			\caption{the distribution of the transmission grid in China. Note: Ref \cite{66}}
			\label{fig:screenshot011}
		\end{minipage}
	\end{figure}

	In China, the development and policymaking in the electricity sector are influenced by a host of complex factors, including government subsidy policies, the competitive structure of the electricity market, consumer demand, and technological advancements, among others \cite{haha1}. The complex interplay of these elements within the electricity market impacts the industry's development. The market's intricate dynamics make predicting its direction difficult, presenting challenges for policy-making and innovation.
	With the continuous development and application of new conductive materials in the future, the game between electricity companies and the government regarding new transmission equipment will gradually surface. Specifically, the theoretical feasibility of superconductivity may accelerate the transformation of the electricity industry, and the emergence of new conductive materials will bring substantial changes and opportunities to the electricity industry \cite{lk}. However, companies and the government often have differing interests and considerations when facing this transformation.
	In the electricity industry, policymakers and leaders must understand and manage the complex interactions and dynamics between the government and companies, which are crucial components of the system.For instance, policymakers need to consider how to allocate and utilize government resources reasonably to encourage and support the development and application of new technologies. Simultaneously, electricity companies need to assess the economic benefits of new technologies and their impacts on existing business models and market shares. Moreover, companies must consider how new technologies will affect their competitive position in the market and how to most effectively use new technologies to enhance their operational efficiency and customer satisfaction.

	The high-quality development of electricity companies depends not only on their own considerations in selecting new equipment but also on government subsidy policies. The government typically supports the promotion and application of new transmission equipment through direct subsidies and tax incentives, among other policy measures. In the electricity industry, government subsidies can impact equipment upgrades. This is because adopting new technologies or replacing equipment may lead to a short-term profit decrease. Thus, when subsidy policy changes, stakeholders scrutinize policy impact on upgrade decisions. The decision-making process then becomes a game between companies and the government. Game theory offers a framework for studying decision-making among actors. Therefore, it is essential to utilize game theory to examine the relationship between power companies and the government and formulate policies before new conductive materials are widely applied.
	
	The study aims to develop government guidelines and promote the construction of a new power system by examining the relationship between government subsidies and electricity companies' equipment upgrade decisions in the event of breakthroughs in superconducting materials technology. In this analytical framework, electricity companies face a strategic decision: whether to upgrade electrical equipment and technology, while the government needs to decide whether to subsidize these upgraded devices. We will systematically analyze the binary relationship subjects in the decision-making process of government subsidy policies and companies adopting new transmission equipment, and their dynamic game process, from the perspective of indirect evolutionary game theory. Considering the influence of bounded rationality and the learning mechanism, we will focus on the decision-making process, striving to solve multiple equilibrium problems. We define an ideal policy event, i.e., the government successfully promoting the widespread application of new electrical equipment, and analyze the various factors that influence the occurrence of this event. We use the analytical tools of indirect evolutionary game theory to predict the potential conflicts of interest between electricity companies and the government when facing new transmission equipment and policies. This theory assumes that participants are rational given their preferences, but their preferences evolve over time.

	In this paper, we will apply indirect evolutionary game theory to study the following scenario: market consumers make choices based on personal utility maximization, but as time progresses, the system state and the utilities of the company and the government will change, requiring the government and electricity companies to dynamically adjust their behavioral strategies. Our game model evaluates the impact of different strategies on both the government and electricity companies, providing corresponding strategies to resolve conflicts of interest and ultimately promote cooperation and development between companies and the government in the field of new power equipment. To help the government and electricity companies respond more effectively to technological and market changes, this study will provide in-depth insights and practical suggestions by answering the following three core questions:
 	\begin{itemize}
		\item In the context of bounded rationality and learning mechanisms, how do electricity companies and the government adjust their strategies to cope with market and technological changes?
		\item How can indirect evolutionary game theory be applied to understand and predict the reactions of electricity companies and the government in the face of new power technologies and policy changes?
		\item In the development process of the electricity industry, how to ensure the maximization of consumer interests while also satisfying the interests of electricity companies and the government?
	\end{itemize}

	 The rest part of this paper is organized as follows. The second section will review relevant research literature, establishing the background and theoretical foundation of the research question.  In Section 3, we will describe the symbols used in the paper and their meanings. The Hotelling Demand Model, which provides the theoretical underpinning for the analysis of the government subsidies and business decision-making process, is shown in section 4. In sections 5 and 6, we will establish the evolutionary game model involving the government and businesses, aiming to resolve conflicts of interest and propose solutions. These sections will also conduct empirical analysis and simulations to deepen our understanding of the impact of government subsidies and business decisions on the electricity industry. This includes the verification of evolutionary stability and the simulation of evolutionary games with time delay. The seventh section will conduct a sensitivity analysis of the initial subsidy values to ensure the reliability and robustness of the research results. The final section will summarize the main findings of the research and discuss their implications for management, providing practical advice and decision-making support for both the government and businesses. The appendix will provide necessary proofs and detailed analysis processes to enhance the credibility of the research and the transparency of the methodology.

	\section{Related Literature}

	\subsection{Evolutionary Game Theory}

	Human decision-making is a multifaceted process involving a large number of factors. Cooperation and environmental change are two pivotal elements in this process. Cooperation can augment the effectiveness and outcomes of decisions, while the continuous change and uncertainty of the environment introduce challenges to decision making. The intricate interaction between these two elements warrants further exploration to understand their profound impacts on human decision-making behavior.
 
	Evolutionary game theory is widely utilized in the research of economics, sociology, and statistical physics to analyze various factors influencing group behavior formation. This theory has generated rich and impactful research findings \cite{5}\cite{6}\cite{7}\cite{8}. Notably, it has achieved significant success in explaining some phenomena in the process of biological evolution. The most classic research results can be traced back to 1973 when Smith and Price first used evolutionary game theory to explain the struggle behavior between animals and first proposed the concept of evolutionarily stable strategies \cite{9}. Building on Smith and Price's work, further developments in the field have expanded our understanding of evolutionary game theory.
	Evolutionary game theory provides a theoretical framework for examining the interaction between individual behavior and fitness. The theory views the evolutionary process as a dynamic process of strategy selection, where an individual's fitness depends not only on the strategy it adopts but also on its interactions with other individuals \cite{9}. Within the application of game theory, the primary objective is to identify a suitable strategy to resolve existing conflicts or to uncover the optimal decision sequence yielding maximum returns.
	However, in the evolutionary process, the competition between different strategies does not occur simultaneously at an exact point in time. Instead, the set of strategies involved in the competition changes over specific time periods. New strategies continuously emerge, compete with existing ones, and potentially replace once dominant strategies. Therefore, a strategy's success should be judged based on the other strategies it encounters in space and time \cite{10}.
 
	In this research, we adopt the evolutionary game model to delve into the strategic decisions of the government and power companies in a complex environment. Evolutionary game theory provides us with a new framework to understand and analyze the decision-making processes of the government and power companies, in which cooperation, environmental changes, and strategy selection and evolution are important aspects of the decision-making process.
	\subsection{Application of Evolutionary Game Theory in Other Fields}
	
	In the study of the game between companies and customers, Li Xiaodong and others \cite{20} used the evolutionary game analysis method to establish an evolutionary game model of the behavior of manufacturing companies and customers. The research found that the degree of servitization of the company and the degree of participation of the customer will affect the knowledge-based interaction between the manufacturing company and the customer.
	Liu \cite{21} constructed an evolutionary game model based on the participation of manufacturers and e-commerce. By solving the game model, they obtained the dynamic evolution rules of manufacturers choosing green innovation and platforms choosing green financing. They believe that companies should not only consider immediate economic benefits but should actively undertake social and environmental responsibilities, improve the company's social image, and obtain considerable additional benefits in the long run.
	In terms of games between companies, Levi \cite{22} used the method of game theory to study two types of joint ventures and proposed that the effective capability of joint ventures is constrained by the most scarce resources.
	Wang and Shi \cite{23} created an industrial pollution evolutionary game model to analyze the game relationship between local governments and companies. Under the static penalty mechanism, the behavior of both the government and the company is uncertain. However, under the dynamic penalty mechanism, the evolutionary path between the government and the company tends to converge to a stable value, that is, the behavior of both parties gradually becomes stable.This implies that the dynamic punishment mechanism can stabilize the game relationship between the government and companies.
	Zhou \cite{24} used an evolutionary game theory analysis method based on system dynamics to study the impact of policy incentives on the development of the electric vehicle industry. The research divided policy incentives into static incentives and dynamic incentives and compared their performance. The results showed that when implementing dynamic incentive policies, the game reached a stable equilibrium.  This means that the electric vehicle industry can develop stably under dynamic incentive policies. Compared with direct subsidy policies, tax incentive policies perform better in stimulating the production of electric vehicles.

	\subsection{Applicability of Evolutionary Game Theory to This Study}

	Game theory is a mathematical framework for analyzing strategic choices (Coninx, 2018) \cite{11}. It models decision paradigms between individuals using the evolutionary dynamics of classical games and characterizes natural selection by the strategy update process. Individuals interact with the game to gain benefits. The higher their benefits, the easier it is for their strategy to be imitated by others. Gene mutation leads to new strategies. Eventually, the entire system is studied, and the statistical laws of group behavior can be obtained through the analysis of the dynamic evolution process of group strategies (Adami, 2016) \cite{12}. However, classical game theory  neglects the dynamic evolution process of individual game behavior, as well as the interaction and adaptive adjustment between individuals. Therefore, evolutionary game theory is used to supplement and expand the research scope of classical game theory. 
	
	EGT is widely used in some concepts, such as Evolutionarily Stable Strategy (ESS) and Replicator Dynamic Equations, to replace classical game theory (Gintis, 2000; Tuyls and Parsons, 2007) \cite{13}\cite{14}. The main purpose of this trend is to reduce the dependence on classical game theory to better understand the evolution process of game behavior. EGT is a descriptive theory, which is very applicable to this study for the following reasons: EGT provides a more detailed and comprehensive analysis framework, can consider the competition between individuals and the evolution of countermeasures, and thus obtain more convincing results. Firstly, classical game theory assumes that the game players are completely rational and have complete information about the environment (Zhang, 2018) \cite{15}. Evolutionary game theory (EGT), on the premise of bounded rationality, improves traditional game theory, assuming that the game players are bounded rational and in an incomplete information environment, viewing the participants' strategy choices as a dynamic adjustment process \cite{16}. However, in some evolutionary games, ESS is not absolutely present (Wang, 2011) \cite{17}. System dynamics provides an effective auxiliary means for the analysis of dynamic evolutionary games under incomplete information, and more intuitively shows the evolutionary trajectory \cite{18}. Secondly, classical game models are difficult to calculate, and they often cannot describe the best way of behavior in terms of optimality and stability, while EGT has more advantages in analyzing strategy change dynamics (Stanford Encyclopedia of Philosophy, 2009) \cite{19}. As mentioned earlier, EGT is very consistent with the dynamic nature of strategic choices of government and power companies.

	\subsection{Innovation Points}
	
	Game theory has been widely used by researchers to study interactions between market participants, deepening the understanding of various market characteristics. 
	The innovative aspects of this study can be outlined as follows:

	\begin{itemize}
		\item Forecasting the Impact of "LK-99" New Material: We study the potential impact of the "LK-99" new material, which could herald a new revolution in material science, on the power equipment industry. Utilizing a proactive approach, we establish a predictive model that provides valuable insights for both government and companies even before the material is commercially available. This approach deviates from traditional reactive models, promoting better strategic planning and preparedness.
		
		\item Incorporating Consumer Influence: Our model introduces a novel element by incorporating consumer influence. The profits of companies using the new ($\pi_{1}$) and traditional power equipment ($\pi_{2}$) are linked with consumer behavior. We provide a theoretical derivation of $\pi_{1}$ and $\pi_{2}$ under consumer influence. This integration presents a more realistic depiction of market dynamics and offers a more comprehensive perspective on the decision-making process of companies.
		
		\item Inclusion of Time Delay: Time delay is integrated into the game between the government and companies, acknowledging that strategies in real-world scenarios may not respond immediately. This inclusion adds an additional layer of complexity to the model and aligns it more closely with real-world situations.
		
		\item Sensitivity Analysis on Subsidy Amounts: We perform a sensitivity analysis on subsidy amounts, revealing that the government's decisions are more susceptible to changes in subsidies compared to those of companies. We explore potential reasons for this differential influence, providing a deeper understanding of how different actors in the system react to changes in financial incentives.
		
	\end{itemize}
	\section{Symbol Explanation}

	\begin{table}[H]
		\centering
		\caption{Nomenclature}
		\begin{tabularx}{\textwidth}{c|X}
			\hline
			\textbf{Parameters} & \textbf{Explanation} \\
			\hline
			$T$ & Consumer's unit mismatch cost \\
			$V_1$ & Consumer's perceived base value of IPTE \\
			$V_2$ & Consumer's perceived base value of TPTE \\
			$q$ & Consumer's preference for energy saving \\
			$eq_1$ & Environmental utility obtained by the consumer through the IPTE \\
			$eq_2$ & Environmental utility obtained by the consumer through the TPTE \\
			$h_1$ & Government's subsidy to consumers for IPTE\\
			$\mu_1$ & Consumer's sensitivity to subsidies \\
			$p_1$ & Consumer's payment price for the IPTE\\
			$p_2$ & Consumer's payment price for the TPTE \\
			$C_1$ & Cost of IPTE \\
			$C_2$ & Cost of TPTE \\
			$P_1$ & Selling price of IPTE \\
			$P_2$ & Selling price of TPTE \\		
			$\pi_1$ & Profit from IPTE \\
			$\pi_2$ & Profit from TPTE \\
			$u_1$ & Government's revenue when the company adopts IPTE without government subsidies\\
			$u_2$ & Additional government revenue when the company adopts IPTE with government subsidies \\
			$u_3$ & Original government revenue when the company adopts TPTE \\
			$t$ & Basic tax paid by the company \\
			$t_1$ & Tax reduction for companies that adopt IPTE \\
			$t_2$ & Penalty tax for companies that adopt TPTE\\		
			$E$ & Revenue function\\		
			$F(y)$ & Government's replication dynamic equation\\		
			$G(x)$ & Power company's replication dynamic equation\\	
			$g(\beta)$ & Coefficient of fiscal subsidy to companies\\		
			$s$ & Initial amount of government subsidy to companies\\		
			$\tau$ & Time delay\\
			\hline		
			\textbf{Variables} & \textbf{Explanation} \\
			\hline
			$x$ & The probability of the government subsidizing IPTE\\		
			$y$ &The probability of a company adopting IPTE\\		
			
			\hline		
			\textbf{Acronyms} & \textbf{Explanation} \\
			\hline
			IPTE&	Innovative power transmission equipment\\	
			TPTE & 	Traditional power transmission equipment\\		
		\end{tabularx}
	\end{table}

	\section{Hotelling Demand Model}
	
	The Hotelling model, an economic model, is used to describe how consumers determine their demand based on their preferences for products and the prices of these products. We apply this model to the power transmission equipment market to analyze consumer choices between new types of transmission equipment and traditional transmission equipment.

	\subsection{Fundamental Assumptions}
	\begin{itemize}
		\item Assume consumers are uniformly distributed along the Hotelling line within the interval [0,1]. Each consumer demands a type of electrical transmission equipment. The electrical transmission equipment market is bifurcated into two categories: novel transmission equipment (located at position 0 on the Hotelling line); and traditional transmission equipment (located at position 1 on the Hotelling line).
		\item It is posited that all consumers have identical incomes and face unit electricity prices of \(p_1\) and \(p_2\) for the two transmission modalities, respectively. This is a simplified assumption intended to render the model more tractable and comprehensible. In reality, consumers' incomes and the prices they encounter may vary, potentially influencing their purchasing decisions.
		\item It's hypothesized that consumers might experience an inconvenience, denoted as \(TX\), when utilizing the novel transmission equipment, while the inconvenience arising from the traditional transmission equipment is \(T(1-X)\). Here, $T$ represents the unit mismatch cost, elucidating the magnitude of the negative utility stemming from the potential impacts on consumers' lives between the novel and traditional transmission methods. We generally assume that \(TX < T(1-X)\), indicating potential instability factors associated with the novel transmission equipment compared to its traditional counterpart.
	\end{itemize}

	\subsection{Consumer Utility Function}
	Let \(V_1\) represent the intrinsic value of the novel transmission equipment, and \(V_2\) signify the intrinsic value of the traditional transmission equipment. Consumers can derive an environmental utility of \(qe_1\) from the novel transmission equipment and \(qe_2\) from the traditional transmission equipment. \(h_1\) denotes the government's subsidy to consumers for the novel transmission equipment, and \(\mu_1\) represents the sensitivity of consumers to this subsidy. \(p_1\) and \(p_2\) are the consumer payment prices for the novel and traditional transmission equipment, respectively. The consumer utility function encapsulates the consumer's preferences towards the novel and traditional transmission equipment, the environmental utility, inconvenience factors, and the response to government subsidies. Equations \eqref{gongshi1} and \eqref{gongshi2} illustrate the utilities derived from using the novel and traditional transmission equipment, respectively.
	
	\begin{equation}
		\label{gongshi1}
		U_1 = V_1 - p_1 + qe_1 - TX + h_1\mu_1
	\end{equation}
	
	\begin{equation}
		\label{gongshi2}
		U_2 = V_2 - p_2 + qe_2 - T(1 - X)
	\end{equation}

	\subsection{Market Share and Company Profit}
	
	When consumers are indifferent between choosing the novel transmission equipment and the traditional transmission equipment, we have \(U_1 = U_2\). This leads to equation \eqref{fc3}:
	
	\begin{equation}
		\label{fc3}
		V_1 - p_1 + qe_1 - TX + h_1\mu_1 = V_2 - p_2 + qe_2 - T(1 - X)
	\end{equation}

	Since the consumer’s perceived base value of IPTE is identical to the consumer’s perceived base value of TPTE, $V1 = V2$.

\begin{figure}[htpb!]
		\centering
		\includegraphics[width=0.5\linewidth]{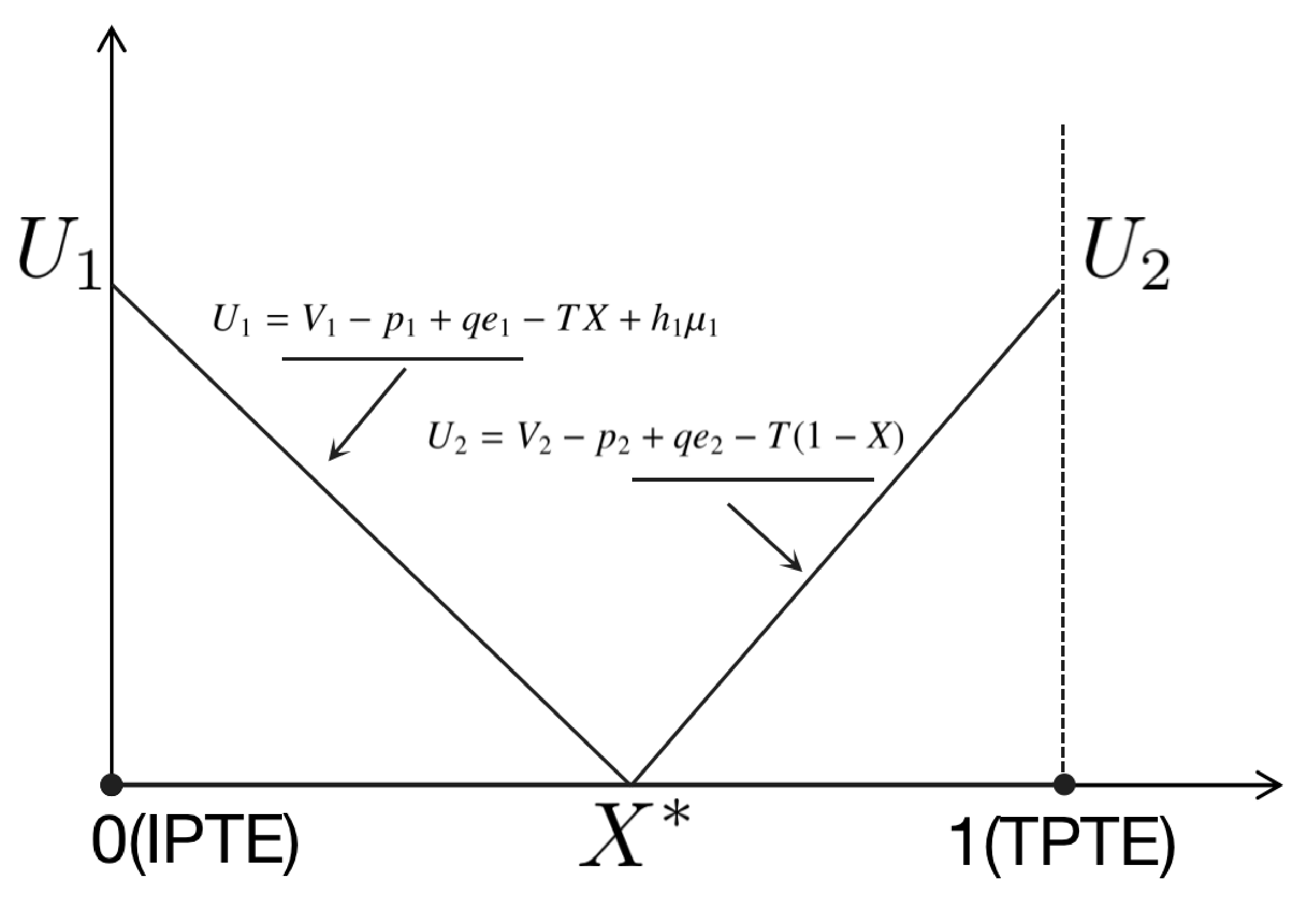}\caption{The market shares of IPTE and TPTE}
		
		\label{061}
	\end{figure}

	From equation \eqref{fc3}, we can deduce the position, denoted as \(X^\ast\), where consumers hold a neutral stance between choosing the novel and traditional transmission equipment:
	
	\begin{equation}
		X^\ast = \frac{[(p2 - p1) + q(e1 - e2) - T + h_1\mu_1]}{2T}
	\end{equation}
    At the point \(X^\ast\), consumers are indifferent between the two options. Thus, we can determine the market shares for the novel and traditional transmission equipment: consumers in the interval [0, \(X^\ast\)] would favor the novel transmission equipment, while those in the interval [\(X^\ast\), 1] would lean towards the traditional transmission equipment, as illustrated in Fig \ref{061}. Subsequently, the profit functions for the companies producing the novel and traditional transmission equipment can be defined as per equations \eqref{gs7} and \eqref{gs8}:
	
	\begin{equation}
		\label{gs7}
		\pi_1 = (P_1-C_1) \times X^\ast
	\end{equation}
	
	\begin{equation}
		\label{gs8}
		\pi_2 = (P_2-C_2) \times (1-X^\ast)
	\end{equation}
	
	Here, \(P_1\) and \(P_2\) respectively represent the prices of the novel and traditional transmission equipment, while \(C_1\) and \(C_2\) are their respective costs.

	\section{Omit this chapter}

			\subsection {Sensitivity Analysis}
			We analyze the evolutionary trends under different initial subsidy amounts $s$ under the conditions of Combination 4 [$u_2 < g(\beta)s + t_1, \pi_2 < \pi_1$].
			
			The parameter settings are as follows:
			$$
			u_2 = 0.5; 
			t_1 = 1; 
			g(\beta) = 1;
			t_2 = 0.5; 
			\pi_1 = 3; 
			\pi_2 = 2; 
			[x_0, y_0] = [0.2, 0.8]
			$$
			
\begin{figure}[H]
				\centering
				\begin{subfigure}{.5\textwidth}
					\centering
					\includegraphics[width=.9\linewidth]{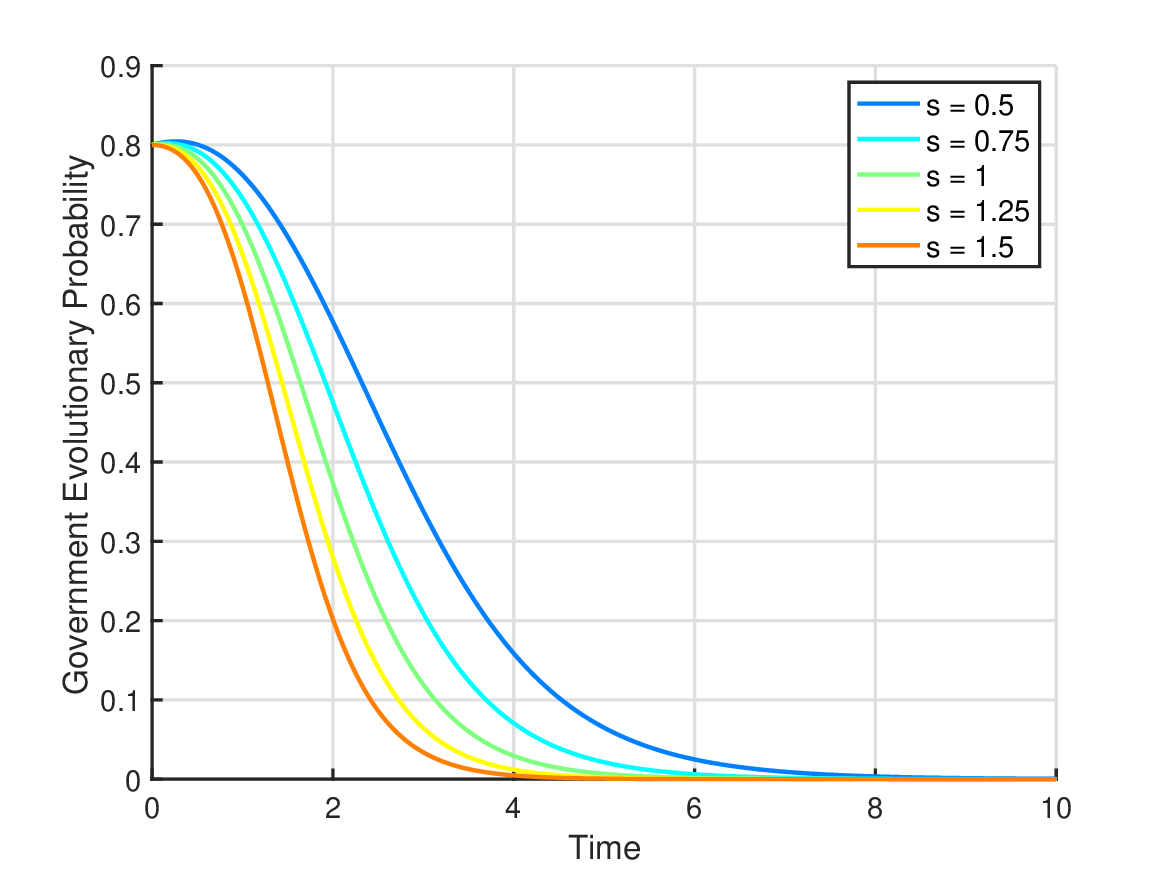}
					\caption{Government evolutionary probability}
					\label{3342}
				\end{subfigure}%
				\begin{subfigure}{.5\textwidth}
					\centering
					\includegraphics[width=.9\linewidth]{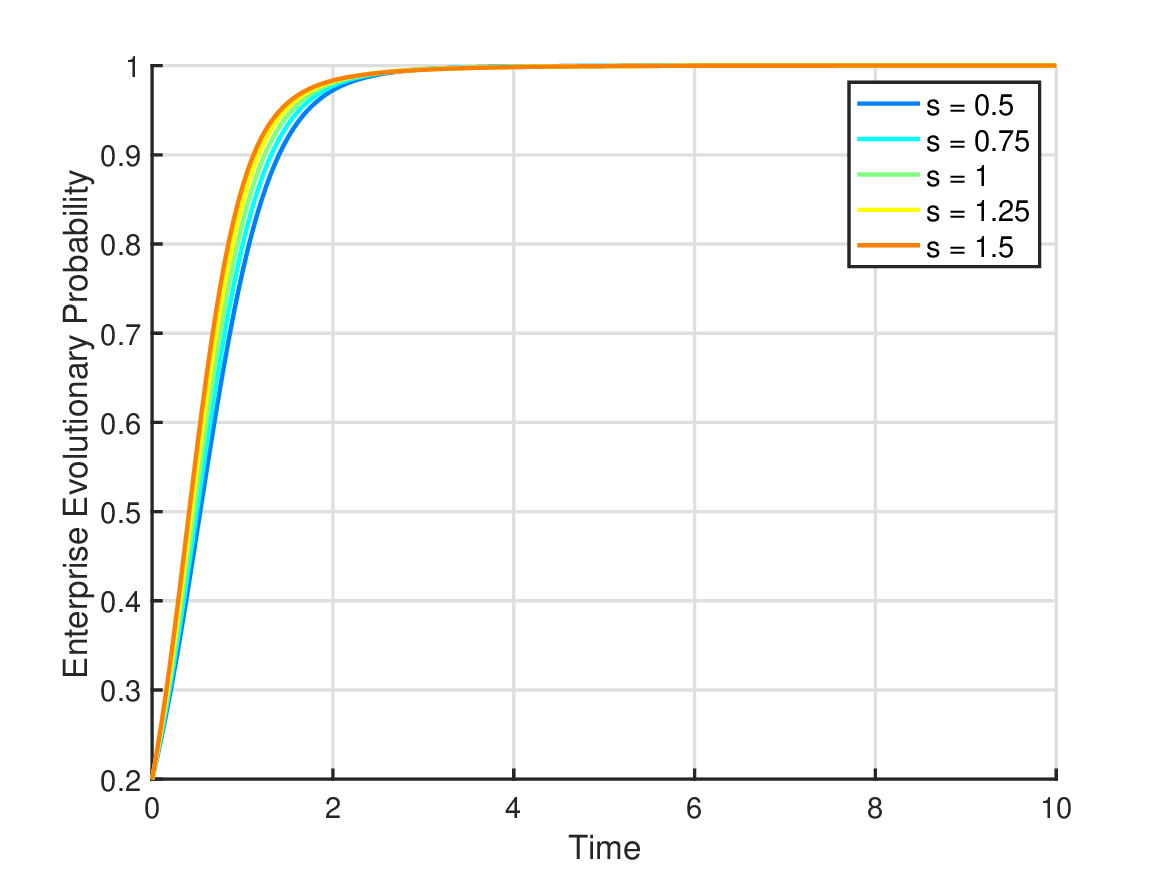}
					\caption{Company  evolutionary probability}
					\label{4244}
				\end{subfigure}
				\caption{Sensitivity analysis}
				\label{fig:subfigures}
			\end{figure}

			We set the initial subsidy amount $s$ to be 0.5, 0.75, 1, 1.25, and 1.5. The evolutionary process under these different $s$ parameters is illustrated in Fig \ref{fig:subfigures}.
			The initial subsidy amount $s$ can considerably affect the dynamics of the game. A lower initial subsidy may lead to slower adoption of new technologies by companies due to the lower initial incentive. In contrast, a higher initial subsidy might lead to a faster adoption of new technologies.
			By conducting numerical simulations under these settings, we can investigate how the initial subsidy amount influences the evolution and stability of the system.

			Through the analysis of the evolutionary process of the government and companies under different initial subsidy amounts, it can be inferred that changes in the initial subsidy amount have a more pronounced effect on the evolution of the government's strategy compared to that of the companies.
			The initial subsidy amount directly influences the government's strategy due to its direct financial implications. An increase in the initial subsidy amount implies a higher initial expenditure for the government, which may lead to a more cautious approach in terms of the government's strategy. On the other hand, a decrease in the initial subsidy could allow the government to be more aggressive in its strategy as it has less financial risk involved.
			Conversely, while companies are also affected by the initial subsidy amount, their strategy evolution is influenced by a combination of factors, including not only the subsidy but also the profitability of adopting new technologies versus maintaining traditional ones. Therefore, the impact of changing the initial subsidy amount on companies may be less direct and less significant compared to that on the government.
			These insights can be instrumental in the formulation of more effective subsidy policies. Policymakers can adjust the initial subsidy amount to influence the strategies of both the government and companies, steering the system towards desirable outcomes.

			\section{Conclusions}

			In the field of material science, the Korean team "LK-99" has developed a potentially transformative material with substantial implications for the power equipment industry. This innovation could fundamentally revolutionize how electricity is generated, transmitted, and consumed. To prepare for this paradigm shift, predictive models are essential to guide government and corporate decision-making.
			
			In answering the first question about bounded rationality and learning mechanisms, our study shows how both the government and electricity companies adjust their strategies over time, particularly when facing market and technological changes. The inclusion of time delays in our model allows us to capture the gradual adjustments that these entities may make.
			For the second question on the application of indirect evolutionary game theory, our research successfully utilizes this framework to predict reactions of electricity companies and the government towards new power technologies and policy changes. Our findings particularly emphasize the impact of initial subsidy amounts, providing a more nuanced understanding of strategic interactions.
			Regarding the third question on balancing consumer interests, while our study has mainly focused on the strategies of government bodies and companies, we have theoretically considered the impact of consumer interests. The current model suggests that the strategies optimized for government and companies also bear the potential to benefit consumers through the faster adoption of more efficient technologies. Future studies should incorporate consumer interests explicitly to create a more comprehensive model.
			
			Due to existing data limitations, the study provides a general overview of anticipated trends. For more precise results, it is recommended to incorporate real-world data into the model once the material exhibits favorable characteristics and begins to be deployed in practical applications.
			
			In conclusion, this research enhances our understanding of the complex interplay among time delay, initial subsidy amount, and strategic choices within the context of an evolutionary game involving government and companies. These findings serve as a valuable theoretical and practical guide for policymakers navigating the complexities of technology adoption and market dynamics. Future research can offer even more nuanced insights by incorporating the role of consumers and adapting to market volatilities.

			\section{Funding}
           Maoqin Yuan was partially supported by Research Foundation of China University of Petroleum-Beijing at Karamay (NO.XQZX20230030) and Talent Project of Tianchi Doctoral Program in Xinjiang Uygur Autonomous Region.
			\section*{Author Contributions}
			Mingyang Li: Conceptualization, methodology, software, validation, formal analysis, investigation and writing (both the original draft and review \& editing).
			
			Maoqin Yuan: Methodology, validation, formal analysis, investigation, writing (both the original draft and review \& editing), supervision, project administration and funding acquisition. In addition,  Maoqin Yuan is the corresponding author, who is responsible for ensuring that the descriptions are accurate and agreed by all authors.
			
			Han Pengsihua: Methodology, investigation, formal analysis, investigation and project administration.
			
			Yuan Yuan: Conceptualization and methodology.
			\nolinenumbers
			\section{Appendix}
			
			Derivation process of the government's replicator dynamic equation \eqref{zf}:
			$$
			\begin{array}{l}
				F(y)=\frac{d y}{d t}\\
				=y\left(E_{21}-E_{2}\right)\\
				=y\{ E_{21}-\left[y E_{21}+(1-y) E_{22}\right]\} \\
				=y\left[(1-y) E_{21}-(1-y) E_{22}\right] \\
				=y(1-y)\left(E_{21}-E_{22}\right)\\
				=y(1-y)[x (u_{2} - t_{1} - g(\beta) s) +  t_{2}(1 - x)]
			\end{array}
			$$
			
			Derivation process of the company's replicator dynamic equation \eqref{gx}:
			$$
			\begin{array}{l}
				G(x)=\frac{d x}{d t}\\
				=x(E_{11}-E_{1}) =x\{E_{11}-\left[x E_{11}+(1-x) E_{12}\right] \}\\
				=x(1-x) E_{11}+(1-x) E_{12} \\
				=x(1-x)\left(E_{11} - E_{12}\right) \\
				=x(1-x)\left\{y\left[g(\beta) s+t_{1}+t_{2}\right]+\left(\pi_{1}-\pi_{2}\right)\right\}
			\end{array}
			$$

\end{document}